\documentclass[11pt]{article}
\usepackage{amsmath}
\usepackage{eucal}
\usepackage{amssymb}

\usepackage{color}
\newtheorem{theorem}{Theorem}[section]

\newtheorem{de}[theorem]{Definition}

\newtheorem{lemma}[theorem]{Lemma}

\begin{document}

\baselineskip 15pt
\title{The influence of $\mathfrak{F_{\mathrm s}}$-quasinormality of subgroups on the structure of finite groups
  \thanks{Research is supported by NNSF of China (Grant \#11071229) and Specialized Research Fund for the Doctoral Program
 of Higher Education of China (Grant \#20113402110036). }}

\author{Xiaolong Yu, Xiaoyu Chen, Wenbin Guo \thanks{Corresponding author}\\
{\small School of Mathematical Sciences, University of Science and Technology of China}\\
{\small Hefei, 230026, P.R. China}\\
{\small Email: wbguo@ustc.edu.cn; yuxiaolong0710@sina.com.cn; jelly@mail.ustc.edu.cn }\\
}

\date{}
\maketitle

\begin{center}
\begin{minipage}{150mm}
{\bf Abstract.} Let $\frak{F}$ be a class
 of finite groups. A subgroup $H$ of a finite group $G$ is said to
 be $\mathfrak{F_{\mathrm s}}$-quasinormal in $G$ if there exists a
normal subgroup $T$ of $G$ such that $HT$ is $s$-permutable in $G$
and $(H\cap T)H_G/H_G$ is contained in the $\frak{F}$-hypercenter
$Z_\infty ^\frak{F} (G/H_G)$ of $G/H_G$. In this paper, we
investigate further the influence of $\mathfrak{F_{\mathrm
s}}$-quasinormality of some subgroups on the structure of finite
groups. New characterization of some classes of finite groups are
obtained.

 \noindent


\end{minipage}
\end{center}

\textbf{This is a revised version of the paper published in Publ. Math. Debrecen (Publ. Math. Debrecen, 82 (3-4): 709-725, 2013). We sincerely thank Professor Luis Ezquerro for pointing out our mistakes.}

\section{Introduction}

Recall that a subgroup $H$ of $G$ is said to be $s$-quasinormal (or
$s$-permutable) in $G$ if $H$ is permutable with every Sylow
subgroup $P$ of $G$ (that is, $HP=PH$). The $s$-permutableity of a
subgroup of a finite group $G$ often yields a wealth of information
about the group $G$ itself. In the past, it has been studied by many
scholars (such as [1-2], [7-9], [13], [17]). Recently, Huang [10]
introduced the following concept:

\begin{de}  Let $\mathfrak{F}$ be a non-empty class of groups and
$H$ a subgroup of a group $G$. $H$ is said to be
$\mathfrak{F_{\mathrm s}}$-quasinormal in $G$ if there exists a
normal subgroup $T$ of $G$ such that $HT$ is $s$-permutable in $G$
and $(H\cap T)H_G/H_G\leq Z_\infty ^\frak{F} (G/H_G)$, where $H_G$
is the maximal normal subgroup of $G$ contained in $H$.
\end{de}

Note that, for a class $\mathfrak{F}$ of groups, a chief factor
$H/K$ of a group $G$ is called $\mathfrak{F}$-central (see [16] or
[4, Definition 2.4.3]) if $[H/K](G/C_G(H/K))\in \mathfrak{F}$. The
symbol $Z_\infty ^\frak{F} (G)$ denotes the
$\mathfrak{F}$-hypercenter of a group $G$, that is, the product of
all such normal subgroups $H$ of $G$ whose $G$-chief factors are
$\mathfrak{F}$-central. A subgroup $H$ of $G$ is said to be
$\mathfrak{F}$-hypercenter in $G$ if $H\leq Z_\infty ^\frak{F}(G)$.

By using this new concept, Huang [10] has given some conditions
under which a finite group belongs to some formations. In this
paper, we will go to further into the influence of
$\mathfrak{F_{\mathrm s}}$-quasinormal subgroups on the structure of
finite groups. New characterizations of some classes of finite
groups are obtained.

All groups considered in the paper are finite and $G$ denotes a
finite group. The notations and terminology in this paper are
standard, as in [4] and [14].

\section{Preliminaries}

Let $\frak{F}$ be a class of finite groups. Then $\frak{F}$ is
called a formation if it is closed under homomorphic image and every
group $G$ has a smallest normal subgroup (called
$\mathfrak{F}$-residual and denoted by $G^\mathfrak{F}$) with
quotient is in $\mathfrak{F}$. $\mathfrak{F}$ is said to be
saturated if it contains every group $G$ with $G/\Phi(G) \in
\mathfrak{F}$. $\mathfrak{F}$ is said to be $S$-closed
($S_n$-closed) if it contains all subgroups (all normal subgroups,
respectively) of all its groups.

We use $\mathfrak{N}$, $\mathfrak{U}$, and $\mathfrak{S}$ to denote
the formations of all nilpotent groups, supersoluble groups and
soluble groups, respectively.

The following known results are useful in our proof.

\textbf{Lemma 2.1} {\rm[8, Lemma 2.2]}. Let $G$ be a group and
$H\leq K\leq G$.

(1) If $H$ is $s$-permutable in $G$, then $H$ is $s$-permutable in
$K$;

(2) Suppose that $H$ is normal in $G$. Then $K/H$ is $s$-permutable
in $G/H$ if and only if $K$ is $s$-permutable in $G$;

(3) If $H$ is $s$-permutable in $G$, then $H$ is subnormal in $G$;

(4) If $H$ and $F$ are $s$-permutable in $G$, the $H\cap F$ is
$s$-permutable in $G$;

(5) If $H$ is $s$-permutable in $G$ and $M\leq G$, then $H\cap M$ is
$s$-permutable in $M$.

\textbf{Lemma 2.2} {\rm[10, Lemma 2.3]}. Let $G$ be a group and
$H\leq K\leq G$.

(1) $H$ is $\mathfrak{F_{\mathrm s}}$-quasinormal in $G$ if and only
if there exists a normal subgroup $T$ of $G$ such that $HT$ is
$s$-permutable in $G$, $H_G\leq T$
 and $H/H_G\cap T/H_G\leq  Z_\infty ^\frak{F} (G/H_G)$;

(2) Suppose that $H$ is normal in $G$. Then $K/H$ is
$\mathfrak{F_{\mathrm s}}$-quasinormal in $G/H$ if and only if $K$
is $\mathfrak{F_{\mathrm s}}$-quasinormal in $G$;

(3) Suppose that $H$ is normal in $G$. Then, for every
$\mathfrak{F_{\mathrm s}}$-quasinormal subgroup $E$ of $G$
satisfying ($|H|$,$|E|$)=1, $HE/H$ is $\mathfrak{F_{\mathrm
s}}$-quasinormal in $G/H$;

(4) If $H$ is $\mathfrak{F_{\mathrm s}}$-quasinormal in $G$ and
$\frak{F}$ is S-closed, then $H$ is $\mathfrak{F_{\mathrm
s}}$-quasinormal in $K$;

(5) If $H$ is $\mathfrak{F_{\mathrm s}}$-quasinormal in $G$, $K$ is
normal in $G$ and $\frak{F}$ is $S_n$-closed, then $H$ is
$\mathfrak{F_{\mathrm s}}$-quasinormal in $K$;

(6) If $G\in \frak{F}$, then every subgroup of $G$ is
$\mathfrak{F_{\mathrm s}}$-quasinormal in $G$.

\textbf{Lemma 2.3} {\rm[6, Lemma 2.2]}. If $H$ is a $p$-subgroup of
$G$ for some prime $p$ and $H$ is $s$-permutable in $G$, then:

(1) $H\leq O_{p}(G)$;

(2) $O^{p}(G)\leq N_{G}(H)$.

\textbf{Lemma 2.4 }{\rm[18]}. If $A$ is a subnormal subgroup of a
group $G$ and $A$ is a $\pi$-group, then $A\leq O_\pi(G)$.

\textbf{Lemma 2.5} {\rm[15, II, Lemma 7.9]}. Let $N$ be a nilpotent
normal subgroup of $G$. If $N\neq 1$ and $N\cap \Phi(G)=1$, then $N$
is a direct product of some minimal normal subgroups of $G$.

\textbf{Lemma 2.6} {\rm[5, Lemma 2.3]}. Let $\mathfrak{F}$ be a
saturated formation containing $\mathfrak{U}$ and $G$ a group with a
normal subgroup $E$ such that $G/E\in \mathfrak{F}$. If $E$ is
cyclic, then $G\in \mathfrak{F}$.

Recall that a subgroup $H$ of $G$ is said to be
$\mathfrak{F}$-supplemented in $G$ if there exists a subgroup $T$ of
$G$ such that $G=HT$ and $T\in \mathfrak{F}$, where $\mathfrak{F}$
is some class of groups. The following Lemma is clear.

\textbf{Lemma 2.7} Let $\mathfrak{F}$ be a formation and $H$ a
subgroup of $G$. If $H$ has an $\mathfrak{F}$-supplement in $G$,
then:

(1) If $N\unlhd G$, then $HN/N$ has an $\mathfrak{F}$-supplement in
$G/N$.

(2) \textbf{If $\mathfrak{F}$ is $S$-closed and $H\leq K\leq G$, then $H$ has an $\mathfrak{F}$-supplement in
$K$.}

\textbf{Lemma 2.8} {\rm[10, Theorem 3.1]}. Let $\mathfrak{F}$ be an
$S$-closed saturated formation containing $\mathfrak{U}$ and $G$ a
group. Then $G\in \mathfrak{F}$ if and only if $G$ has a normal
subgroup $E$ such that $G/E \in \mathfrak{F}$ and every maximal
subgroup of every non-cyclic Sylow subgroup of $E$ not having a
supersoluble supplement in $G$ is $\mathfrak{U}_{\mathrm
s}$-quasinormal in $G$.

\textbf{Lemma 2.9} {\rm[10, Theorem 3.2]}. Let $\mathfrak{F}$ be a
saturated formation containing $\mathfrak{U}$ and $G$ a group. Then
$G\in \mathfrak{F}$ if and only if $G$ has a soluble normal subgroup
$E$ such that $G/E \in \mathfrak{F}$ and every maximal subgroup of
every non-cyclic Sylow subgroup of $F(E)$ not having a supersoluble
supplement in $G$ is $\mathfrak{U}_{\mathrm s}$-quasinormal in $G$.

\textbf{Lemma 2.10} {\rm[3, Main Theorem]}. Suppose $G$ has a Hall
$\pi$-subgroup and $2\notin \pi$. Then all the Hall $\pi$-subgroups
are conjugate in $G$.

\textbf{Lemma 2.11} {\rm[6, Lemma 2.5]}. Let $G$ be a group and $p$
a prime such that $p^{n+1}\nmid |G|$ for some integer $n\geq 1$. If
$(|G|,(p-1)(p^{2}-1)\cdot \cdot \cdot (p^{n}-1))=1$, then $G$ is
$p$-nilpotent.

The generalized Fitting subgroup $F^{*}(G)$ of a group $G$ is the
product of all normal quasinilpotent subgroups of $G$. We also need
in our proofs the following well-known facts about this subgroups
(see [12, Chapter X]).

\textbf{Lemma 2.12}. Let $G$ be a group and $N$ a subgroup of $G$.

(1) If $N$ is normal in $G$, then $F^{*}(N)\leq F^{*}(G)$.

(2) If $N$ is normal in $G$ and $N\leq F^{*}(G)$, then
$F^{*}(G)/N\leq F^{*}(G/N)$.

(3) $F(G)\leq F^{*}(G)=F^{*}(F^{*}(G))$. If $F^{*}(G)$ is soluble,
then $F^{*}(G)=F(G)$.

(4) $C_{G}(F^{*}(G))\leq F(G)$.

(5) $F^{*}(G)=F(G)E(G)$, $F(G)\cap E(G)=Z(E(G))$ and $E(G)/Z(E(G))$
is the direct product of simple non-abelian groups, where $E(G)$ is
the layer of $G$.

\textbf{Lemma 2.13} {\rm[8, Lemma 2.15-2.16]}.

(1) If $H$ is a normal soluble subgroup of a group $G$, then
$F^{*}(G/\Phi (H))=F^{*}(G)/\Phi (H)$.

(2) If $K$ is a normal $p$-subgroup of a group $G$ contained in
$Z(G)$, then $F^{*}(G/K)=F^{*}(G)/K$.

\section{New Characterization of supersoluble groups}

\begin{lemma} Let $\textit{p}$ be the smallest prime dividing $|G|$
and $P$ some Sylow $\textit{p}$-subgroup of $G$. Then $G$ is soluble
if and only if every maximal subgroup of $P$ is
$\mathfrak{S_{\mathrm s}}$-quasinormal in $G$.
\end{lemma}

{\noindent\bf Proof.} The necessity is obvious since $Z_\infty
^\frak{S}(G)=G$ whenever $G\in \mathfrak{S}$. Hence we only need to
prove the sufficiency. Suppose that the assertion is false and let
$G$ be a counterexample of minimal order. Then $\textit{p}=2$ by the
well known Feit-Thompson Theorem of groups of odd order. We proceed
the proof via the following steps:

\textsl{(1) $O_{2}(G)=1$.}

Assume that $N=O_{2}(G)\neq 1$. Then $P/N$ is a Sylow 2-subgroup of
$G/N$. Let $M/N$ be a maximal subgroup of $P/N$. Then $M$ is a
maximal subgroup of $P$. By the hypothesis and Lemma 2.2(2), $M/N$
is $\mathfrak{S_{\mathrm s}}$-quasinormal in $G/N$. The minimal
choice of $G$ implies that $G/N$ is soluble. It follows that $G$ is
soluble, a contradiction. Hence (1) holds.

\textsl{(2) $O_{2^{'}}(G)=1$.}

Assume that $D=O_{2^{'}}(G)\neq 1$. Then $PD/D$ is a Sylow
2-subgroup of $G/D$. Suppose that $M/D$ is a maximal subgroup of
$PD/D$. Then there exists a maximal subgroup $P_{1}$ of $P$ such
that $M=P_{1}D$. By the hypothesis and Lemma 2.2(3), $M/D=P_{1}D/D$
is $\mathfrak{S_{\mathrm s}}$-quasinormal in $G/D$. Hence $G/D$ is
soluble by the choice of $G$. It follows that $G$ is soluble, a
contradiction.

\textsl{(3) Final contradiction.}

Let $P_{1}$ be a maximal subgroup of $P$. By the hypothesis, there
exists a normal subroup $K$ of $G$ such that $P_{1}K$ is
$s$-permutable in $G$ and $(P_{1}\cap K)(P_{1})_{G}/(P_{1})_{G}\leq
Z_\infty ^\frak{S}(G/(P_{1})_{G})$. Note that $Z_\infty
^\frak{S}(G)$ is a soluble normal subgroup of $G$. By (1) and (2),
we have $(P_{1})_{G}=1$ and $Z_\infty ^\frak{S}(G)=1$. This induces
that $P_{1}\cap K=1$. If $K=1$, then $P_{1}$ is $s$-permutable in
$G$ and so $P_{1}=1$ by (1) (2) and Lemma 2.3(1). This means that
$|P|=2$. Then by [14, (10.1.9)], $G$ is $2$-nilpotent and so $G$ is
soluble, a contradiction. We may, therefore, assume that $K\neq 1$.
If $2\mid |K|$, then $|K_{2}|=2$, where $K_{2}$ is a Sylow
2-subgroup of $K$. By [14, (10.1.9)] again, we see that $K$ is
2-nilpotent, and so $K$ has a normal 2-complement $K_{2^{'}}$. Since
$K_{2^{'}}$ char $K\unlhd G$, $K_{2^{'}}\unlhd G$. Hence by (2),
$K_{2^{'}}=1$. Consequently $|K|=2$, which contradicts (1). If
$2\nmid |K|$, then $K$ is a $2^{'}$-group. Hence by (2), $K\leq
O_{2^{'}}(G)=1$, also a contradiction. This completes the proof.

\begin{theorem}Let $G=AB$, where $A$ is a subnormal subgroup of
$G$, and $B$ is a supersoluble Hall subgroup of $G$ in which all
Sylow subgroups are cyclic. If every maximal subgroup of every
non-cyclic Sylow subgroup of $A$ is $\mathfrak{U_{\mathrm
s}}$-quasinormal in $G$, then $G$ is supersoluble.
\end{theorem}

{\noindent\bf Proof.}  Suppose that the assertion is false and let
$G$ be a counterexample of minimal order. Then:

\textsl{(1) Each proper subgroup of $G$ containing $A$ is
supersoluble.}

Let $A\leq M< G$. Then $M=M\cap AB=A(M\cap B)$. Obviously, $M\cap B$
is a Hall subgroup of $M$ and every Sylow subgroup of $M\cap B$ is
cyclic. By Lemma 2.2(4), every maximal subgroup of every non-cyclic
Sylow subgroup of $A$ is $\mathfrak{U_{\mathrm s}}$-quasinormal in
$M$. The minimal choice of $G$ implies that $M$ is supersoluble.

\textsl{(2) Let $H$ be a non-trivial normal $p$-subgroup of $G$ for
some prime $p$. If $H$ contains some Sylow $p$-subgroup of $A$ or a
Sylow $p$-subgroup of $A$ is cyclic or $H\leq A$, then $G/H$ is
supersoluble.}

If $A\leq H$, then $G/H=BH/H\cong B/(B\cap H)$ is supersoluble. Now
we can assume that $A\nleq H$. Clearly, $G/H=(AH/H)(BH/H)$, where
$AH/H$ is subnormal in $G/H$ and $BH/H$ is supersoluble. Let $Q/H$
be any non-cyclic Sylow $q$-subgroup of $AH/H$ and $Q_{1}/H$ a
maximal subgroup of $Q/H$. Then there exists a non-cyclic Sylow
$q$-subgroup $A_{q}$ of $A$ such that $Q=A_{q}H$ and a maximal
subgroup $A_{1}$ of $A_{q}$ such that $Q_{1}=A_{1}H$. If $H\leq A$,
then the assertion holds by the choice of $G$ and Lemma 2.2(2). We
may, therefore, assume that $H\nleq A$. Let $P$ be a Sylow
$p$-subgroup of $A$. Assume that $P$ is cyclic or $P\leq H$. Then
$p\neq q$. Clearly, $Q_{1}\cap A_{q}=A_{1}$ is a maximal subgroup of
$A_{q}$. By the hypothesis, $A_{1}$ is $\mathfrak{U_{\mathrm
s}}$-quasinormal in $G$. Therefore, $Q_{1}/H=A_{1}H/H$ is
$\mathfrak{U_{\mathrm s}}$-quasinormal in $G/H$ by Lemma 2.2(3).
This shows that the conditions of the theorem are true for $G/H$ and
so $G/H$ is supersoluble by the minimal choice of $G$.

\textsl{(3) There exists at least one Sylow subgroup of $A$ which is
non-cyclic.}

It follows from the well known fact that a group $G$ is supersoluble
if all its Sylow subgroups are cyclic.

\textsl{(4) $G$ is soluble.}

If $A\neq G$, then $A$ is supersoluble by (1). Let $p$ be the
largest prime divisor of $|A|$. Then $A_{p}\unlhd A$. By Lemma 2.4,
$A_{p}\leq O_{p}(G)$. By (2), $G/O_{p}(G)$ is supersoluble. It
follows that $G$ is soluble.

We now only need to consider the case that $A=G$. If $G$ is not
soluble and let $p$ be the minimal prime divisor of $|G|$. Then
$p=2$ by the well-known Feit-Thompson Theorem. Hence by Lemma 3.1,
$G$ is soluble.

\textsl{(5) $G$ has a unique minimal normal subgroup $N$ such that
$N=O_p(G)=C_G(N)$ is a non-cyclic $p$-subgroup of $G$ for some prime
$p$ and $G=[N]M$, where $M$ is a supersoluble maximal subgroup of
$G$.}

Let $N$ be an arbitrary minimal normal subgroup of $G$. By (4), $N$
is a $p$-group. If $p\in \pi (B)$, then the Sylow $p$-subgroups of
$G$ are cyclic and so the Sylow $p$-subgroups of $A$ are cyclic. If
$p\notin \pi (B)$, then clearly, $N\subseteq A$. Hence by (2), $G/N$
is supersoluble. If $N$ is cyclic, then by Lemma 2.6, $G$ is
supersoluble, a contradiction. Since the class of all supersoluble
groups is a saturated formation, $N$ is the only minimal normal
subgroup $N$ of $G$ and $\Phi(G)=1$. This implies that (5) holds.

\textsl{(6) $N$ is not a Sylow subgroup of $G$ and $Z_\infty
^\mathfrak{U}(G)=1$.}

By (5), clearly, $Z_\infty ^\mathfrak{U}(G)=1$. Assume that $N$ is a
Sylow $p$-subgroup of $G$. Let $N_{1}$ be a maximal subgroup of $N$.
Then by hypothesis, $N_{1}$ is $\mathfrak{U_{\mathrm
s}}$-quasinormal in $G$. Hence there exists a normal subgroup $K$ of
$G$ such that $N_{1}K$ is $s$-permutable in $G$ and $N_{1}\cap K\leq
Z_\infty ^\frak{U}(G)=1$ since $(N_{1})_{G}=1$. It follows that
$N_{1}\leq N_{1}\cap N\leq N_{1}\cap K=1$. Hence $|N|=p$. This
contradiction shows that $N$ is not a Sylow $p$-subgroup of $G$.

\textsl{(7) $A$ is supersoluble.}

If $A$ is not supersoluble, then $G=A$ by (1). Let $q$ be the
largest prime divisor of $|G|$ and $Q$ is a Sylow $q$-subgroup of
$G$. Then $QN/N$ is a Sylow $q$-subgroup of $G/N$. Since $G/N$ is
supersoluble, $QN/N\unlhd G/N$. It follows that $QN\unlhd G$. Let
$P$ be a non-cyclic Sylow $p$-subgroup of $G=A$. If $p=q$, then
$P=Q=QN\unlhd G$. Therefore $N=O_p(G)=P$ is the Sylow $p$-subgroup
of $G$, a contradiction. Assume that $q> p$. Then clearly $QP=QNP$
is a subgroup of $G$. Since $N\nleq \Phi (G)$, $N\nleq \Phi (P)$ by
[11, III, Lemma 3.3(a)]. Let $P_1$ be a maximal subgroup of $P$ such
that $N\nleq P_{1}$.  Then $(P_1)_G=1$. By the hypothesis, $P_1$ is
$\mathfrak{U_{\mathrm s}}$-quasinormal in $G$. Hence, there exists a
normal subgroup $T$ of $G$ such that $P_1T$ is $s$-permutable in $G$
and $P_1\cap T\leq Z_\infty^\mathfrak{U}(G)=1$. Obviously, $T\neq 1$
(In fact, if $T=1$, then $P_{1}\leq O_{p}(G)=N$ by Lemma 2.3(1).
Hence $P_{1}=N$ or $P=N$. This is impossible). Thus $N\leq T$, and
so $P_1\cap N\leq P_1\cap T=1$. This induces that $|N|=|P:P_{1}|=p$,
which contradicts (5). Thus (7) holds.

\textsl{(8) The final contradiction.}

Let $q$ be the largest prime divisor of $|A|$ and $A_{q}$ a Sylow
$p$-subgroup of $A$. Since $A$ is supersoluble by (7), $A_{q}\unlhd
A$. Hence $A_{q}\leq O_{q}(G)$. If $q\mid |B|$, then $O_{q}(G)\leq
G_{q}$, where $G_{q}$ is a cyclic Sylow $q$-subgroup of $B$ and so
$O_{q}(G)$ is cyclic. In view of (2), $G/O_{q}(G)$ is supersoluble.
It follows that $G$ is supersoluble, a contradiction. Hence $q\nmid
|B|$. Then, $A_{q}$ is a Sylow $q$-subgroup of $G$ and so
$A_{q}=O_{q}(G)\neq 1$. This means that $q=p$ and so
$N=A_{p}=G_{p}$, which contradicts (6). The final contradiction
completes the proof.

\begin{theorem}Let $\mathfrak{F}$ be an $S$-closed saturated formation
containing $\mathfrak{U}$ and $H$ a normal subgroup of $G$ such that
$G/H \in \mathfrak{F}$. Suppose that every maximal subgroup of every
non-cyclic Sylow subgroup of $F^{*}(H)$ having no supersoluble
supplement in $G$ is $\mathfrak{U_{\mathrm s}}$-quasinormal in
 $G$. Then $G\in \mathfrak{F}$.
\end{theorem}

{\noindent\bf Proof.} We first prove that the theorem is true if
$\mathfrak{F}=\mathfrak{U}$. Suppose that the assertion is false and
consider a counterexample for which $|G||H|$ is minimal. Then:

  \textsl{(1) $H=G$ and $F^{*}(G)=F(G)$.}

  By Lemma 2.8, $F^{*}(H)$ is supersoluble. Hence $F^{*}(H)=F(H)$ by Lemma
  2.12(3). Since $(H,H)$ satisfies the hypothesis, the minimal choice of
  $(G,H)$ implies that $H$ is supersoluble if $H< G$. Then $G\in
  \mathfrak{U}$ by Lemma 2.9, a contradiction.

  \textsl{(2) Every proper normal subgroup $N$ of $G$ containing
  $F^{*}(G)$ is supersoluble.}

  Let $N$ be a proper normal subgroup of $G$ containing $F^{*}(G)$.
  By Lemma 2.12, $F^{*}(G)=F^{*}(F^{*}(G))\leq F^{*}(N)\leq
  F^{*}(G)$. Hence $F^{*}(N)=F^{*}(G)$. Let $M$ be a maximal
  subgroup of any non-cyclic Sylow subgroup of $F^{*}(N)$. If there
  exists a supersoluble subgroup $T$ such that $G=MT$, then
  $N=M(N\cap T)$ and $N\cap T\in \mathfrak{U}$. This means that $M$
  has a supersoluble supplement in $N$. Now assume that $M$ has no
  supersoluble supplement in $G$. Then by hypothesis and Lemma
  2.2(4), $M$ is $\mathfrak{U_{\mathrm s}}$-quasinormal in $N$. This shows that
  $(N,N)$ satisfies the hypothesis. Hence $N$ is supersoluble by the
  minimal choice of $(G,H)$.

  \textsl{(3) If $p\in \pi(F(G))$, then $\Phi(O_{p}(G))=1$ and so $O_{p}(G)$
  is elementary abelian. In particular, $F^{*}(G)=F(G)$ is abelian
  and $C_{G}(F(G))=F(G)$.}

  Suppose that $\Phi(O_{p}(G))\neq 1$ for some $p\in \pi(F(G))$. By
  Lemma 2.13(1), we have $F^{*}(G/\Phi (O_{p}(G)))=F^{*}(G)/\Phi
  (O_{p}(G))$. By using Lemma 2.2, we see that the pair $(G/\Phi (O_{p}$
  $(G)),F^{*}(G)/\Phi (O_{p}(G)))$ satisfies the hypothesis. The
  minimal choice of $(G,H)$ implies $G/\Phi (O_{p}(G))$
  $\in \mathfrak{U}$. Since $\mathfrak{U}$ is a saturated formation, we
  obtain that $G\in \mathfrak{U}$, a contradiction. This means that
  $\Phi(O_{p}(G))=1$ and so $O_{p}(G)$ is elementary abelian. Hence $F^{*}(G)=F(G)$
  is abelian and $F(G)\leq C_{G}(F(G))$. Put $N=C_{G}(F(G))$. Then, clearly, $F(N)=F(G)$. If
  $N=G$, then $F(G)\leq Z(G)$. Let $P_1$ be a maximal subgroup of some Sylow $p$-subgroup
  of $F(G)$. Then $F(G/P_1)=F(G)/P_1$ by Lemma 2.13(2). Hence $(G/P_1,
  F(G)/P_1)$ satisfies the hypothesis and so $G/P_1\in \mathfrak{F}$.
  Then since $P\leq Z(G)$, we obtain $G\in \mathfrak{F}$. This contradiction shows that
  $N<G$. Hence by (2), $N$ is soluble and so $C_{N}(F(N))\subseteq F(N)$.
  It follows that $N=C_{G}(F(G))=F(G)$.

 \textsl{ (4) $G$ has no normal subgroup of prime order contained in
$F(G)$.}

Suppose that $L$ is a normal subgroup of $G$ contained in $F(G)$ and
$|L|=p$. Put $C=C_{G}(L)$. Clearly, $F(G)\leq C\unlhd G$. If $C<G$,
then $C$ is soluble by (2). Since $G/C$ is cyclic, $G$ is soluble.
Then by the hypothesis and Lemma 2.9, $G\in \mathfrak{U}$, a
contradiction. Hence $C=G$ and so $L\leq Z(G)$. By Lemma 2.13(2)
$F^{*}(G/L)=F^{*}(G)/L=F(G)/L$. Hence $G/L$ satisfies the hypothesis
by Lemma 2.2. The minimal choice of $(G,H)$ implies that $G/L\in
\mathfrak{U}$ and consequently $G$ is supersoluble, a contradiction.

\textsl{(5) For some $p\in \pi(F(G))$, $O_{p}(G)$ is a non-cyclic
Sylow $p$-subgroup of $F(G)$.}

Clearly, $F(G)=O_{p_{1}}(G)\times O_{p_{2}}(G)\times \cdot \cdot
\cdot \times O_{p_r}(G)$ for some primes $p_i$, $i=1,2,\cdots , r.$
If all Sylow subgroups of $F(G)$ are cyclic, then
$G/C_{G}(O_{p_{i}}(G))$ is abelian for any $i\in \{1\cdot \cdot
\cdot r\}$ and so
$G/\cap^{r}_{i=1}C_{G}(O_{p_{i}}(G))=G/C_{G}(F(G))=G/F(G)$ is
abelian. Therefore $G$ is soluble. It follows from Lemma 2.9 and the
hypothesis that $G\in \mathfrak{U}$, a contradiction.

\textsl{(6) Every maximal subgroup of every non-cyclic Sylow
subgroup of $F(G)$ has no supersoluble supplement in $G$.}

Let $P$ be a non-cyclic Sylow subgroup of $F(G)$ and $P_{1}$ a
maximal subgroup of $P$. Then $P=O_{p}(G)$ for some $p\in
\pi(F(G))$. If $P_{1}$ has a supersoluble supplement in $G$, that
is, there exists a supersoluble subgroup $K$ of $G$ such that
$G=P_{1}K=O_{p}(G)K$, then $G/O_{p}(G)\simeq K/K\cap O_{p}(G)$ is
supersoluble and so $G$ is soluble. Hence as above, $G\in
\mathfrak{U}$, a contradiction.

\textsl{(7) $P\cap \Phi(G)\neq 1$, for some non-cyclic Sylow
subgroup $P$ of $F(G)$.}

Assume that $P\cap \Phi(G)=1$. Then $P=R_{1}\times R_{2}\times$
$\cdot \cdot \cdot \times R_{m}$, where $R_{i}(i\in \{1,\cdot \cdot
\cdot m\})$ is a minimal normal subgroup of $G$ by Lemma 2.5. We
claim that $R_{i}$ are of order $\textit{p}$ for all $i\in \{1,\cdot
\cdot \cdot m\}$. Assume that $|R_{i}|>\textit{p}$, for some
$\textit{i}$. Without loss of generality, we let
$|R_{1}|>\textit{p}$. Let $R^{*}_{1}$ be a maximal subgroup of
$R_{1}$. Obviously, $R^{*}_{1}\neq 1$. Then $R^{*}_{1}\times
R_{2}\times$ $\cdot \cdot \cdot \times
   R_{m}=P_{1}$ is a maximal subgroup of $P$. Put $T=R_{2}\times$
   $\cdot \cdot \cdot \times
   R_{m}$, Clearly $(P_{1})_{G}=T$. By (6) and the hyperthesis,
   $P_{1}$ is $\mathfrak{U_{\mathrm s}}$-quasinormal in $G$. Hence by Lemma
   2.2(1), there exists a normal subgroup $N$ of $G$ such that $(P_{1})_{G}\leq
   N$, $P_{1}N$ is $s$-permutable in $G$ and $P_{1}/(P_{1})_{G} \cap N/(P_{1})_{G}\leq Z_\infty
^\frak{U} (G/(P_{1})_{G})$. Assume that $P_{1}/(P_{1})_{G} \cap
N/(P_{1})_{G}\neq 1$. Let $Z_\infty ^\frak{U}
(G/(P_{1})_{G})=V/(P_{1})_{G}=V/T$. Then $P/T \cap V/T \unlhd G/T$.
Since $P\cap V\geq P_{1}\cap N\cap V\geq P_{1}\cap N>
(P_{1})_{G}=T$, we have $P/T \cap V/T\neq 1$. Because $P/T\simeq
R_{1}$ and $R_{1}$ is a minimal normal subgroup of $G$,
$P/T\subseteq V/T$. This implies that $|R_{1}|=|P/T|=p$. This
contradiction shows that $P_{1}\cap N=(P_{1})_{G}=T$. Consequently
$P_{1}N=R^{*}_{1}TN=R^{*}_{1}N$ and $R^{*}_{1}\cap N=1$. Since
$R_{1}\cap N\unlhd G$, $R_{1}\cap N=1$ or $R_{1}\cap N=R_{1}$. But
since $R^{*}_{1}\cap N=1$, we have that $R_{1}\cap N=1$. Thus
$R^{*}_{1}=R^{*}_{1}(R_{1}\cap N)=R_{1}\cap R^{*}_{1}N$ is
$s$-permutable in $G$. It follows from Lemma 2.3(2) that
$O^{p}(G)\leq N_{G}(R_{1}^{*})$. Thus $|G:N_{G}(R_{1}^{*})|$ is a
power of $p$ for every maximal subgroup $R_{1}^{*}$ of $R_{1}$. This
induces that $p$ divides the number of all maximal subgroups of
$R_{1}$. This contradicts [11, III, Theorem 8.5(d)]. Therefore
$|R_{i}|=p$, which contradicts (4). Thus (7) holds.

\textsl{(8) $F(G)=P$ is a $p$-group, $P$ contains a unique minimal
normal subgroup $L$ of $G$ and $L\subseteq \Phi (G)$.}

Suppose that $1\neq Q$ is a Sylow $q$-subgroup of $F(G)$ for some
prime $q\neq p$ and let $L$ be a minimal normal subgroup of $G$
contained in $P\cap \Phi (G)$. By (3), $Q$ is elementary abelian. By
Lemma 2.12, $F^{*}(G/L)=F(G/L)E(G/L)$ and $[F(G/L),E(G/L)]=1$, where
$E(G/L)$ is the layer of $G/L$. Since $L\leq \Phi (G)$,
$F(G/L)=F(G)/L$. Now let $E/L=E(G/L)$. Since $Q$ is normal in $G$
and $[F(G)/L,E/L]=1$, we have $[Q,E]\leq Q\cap L=1$. It follows from
(3) that $F(G)E\leq C_{G}(Q)\unlhd G$. If $C_{G}(Q)< G$, then
$C_{G}(Q)$ is supersoluble by (1) and (2). Thus $E(G/L)=E/L$ is
supersoluble and consequently $F^{*}(G/L)=F(G)/L$ by Lemma 2.12(5).
Now, by Lemma 2.2, we see that $(G/L,F(G)/L)$ satisfies the
hypothesis. The minimal choice of $(G,H)$ implies that $G/L$ is
supersoluble and so is $G$. This contradiction shows that
$C_{G}(Q)=G$, i.e. $Q\leq Z(G)$, which contradicts (4). Thus
$F(G)=P$.

Let $X$ be a minimal normal subgroup of $G$ contained in $P$ with
$X\neq L$. Let $E/L=E(G/L)$ is the layer of $G/L$. As above, we see
that  $F^{*}(G/L)=F(G/L)E(G/L)$ and $[F(G)/L,E/L]=1$. Hence
$[X,E]\leq X\cap L=1$, i.e., $[X,E]=1$. It follows from (3) that
$F(G)E\leq C_{G}(X)\unlhd G$. If $C_{G}(X)<G$, then $C_{G}(X)$ is
supersoluble by (1) and (2). Thus $E(G/L)=E/L$ is supersoluble and
consequently $F^{*}(G/L)=F(G)/L$. Obviously, $G/L$ satisfies the
hypothesis. By the choice of $(G,H)$, we have that $G/L$ is
supersoluble and so is $G$, a contradiction. Hence $C_{G}(X)=G$,
i.e. $X\leq Z(G)$, which also contradicts (4). Thus $L$ is the
unique minimal normal subgroup of $G$ contained in $P$. Finally,
$L\subseteq \Phi (G)$ by (7).

\textsl{(9) $L<P$.}

Suppose $L=P$. Let $P_{1}$ be a maximal subgroup of $P$ such that
$P_1$ is normal in some Sylow subgroup of $G$. Then $(P_{1})_{G}=1$.
By the hypothesis and (8), $P_{1}$ is $\mathfrak{U_{\mathrm
s}}$-quasinormal in $G$. Hence there exists a normal subgroup $K$ of
$G$ such that $P_{1}K$ is $s$-permutable in $G$ and $P_{1}\cap K\leq
Z_\infty ^\mathfrak{U}(G)$. If $P_{1}\cap K\neq 1$, then
$1<P_{1}\cap K\leq P\cap Z_\infty ^\mathfrak{U}(G)$, which implies
that $P=P\cap Z_\infty ^\mathfrak{U}(G)$ and $|P|=p$ since $P$ is a
minimal normal subgroup of $G$. This contradicts (4). So we may
assume $P_{1}\cap K=1$. Since $P$ is a minimal normal subgroup of
$G$, $P\cap K=P$ or $1$. If $P\cap K=P$, then $P\subseteq K$, and so
$|P|=p$, which contradicts (4). If $P\cap K=1$, then $P\cap
P_{1}K=P_{1}(P\cap K)=P_{1}$. Hence $P_{1}$ is $s$-permutable in
$G$. Then by Lemma 2.3(2), $O^{p}(G)\leq N_{G}(P_{1})$. This induces
that $P_{1}\unlhd G$. This means that $P_{1}=(P_{1})_{G}=1$ and
$|P|=p$, also a contradiction.

\textsl{ (10) Final contradiction (for
$\mathfrak{F}=\mathfrak{U}$).}

By (3) and (8), $P$ is an elementary abelian group, and so $L$ has a
complement in $P$, $T$ say. Let $P_{1}=TL_{1}$, where $L_{1}$ is a
maximal subgroup of $L$. Then $1\neq P_{1}$ and clearly $P_1$ is a
maximal subgroup of $P$ such that $P_1$ is  normal in some Sylow
subgroup of $G$. Hence by (6), $P_{1}$ is $\mathfrak{U_{\mathrm
s}}$-quasinormal in $G$ and $(P_{1})_{G}=1$ since $L$ is the unique
minimal normal subgroup of $G$ contained in $P$. Hence there exists
a normal subgroup $S$ of $G$ such that $P_{1}S$ is $s$-permutable in
$G$ and $P_{1}\cap S\leq Z_\infty ^\mathfrak{U}(G)$. If $P_{1}\cap
S\neq 1$, then $1< P_{1}\cap S\leq P\cap Z_\infty ^\mathfrak{U}(G)$
and so $G$ has a minimal normal subgroup $N$ of order $p$ contained
in $P$, which is contrary to (4). Hence $P_{1}\cap S=1$. If $P\cap
S\neq 1$, then $L\leq P \cap S$ and so $L_{1}\leq S$, which
contradicts $P_{1}\cap S= 1$. If $P\cap S=1$, then
$P_{1}=P_{1}(P\cap S)=P\cap P_{1}S$ is $s$-permutable in $G$. Hence
$O^{p}(G)\leq N_{G}(P_{1})$ by Lemma 2.3(2). It follows that
$P_{1}\unlhd G$, which contradicts $(P_{1})_{G}=1$. The final
contradiction shows that the theorem holds when
$\mathfrak{F}=\mathfrak{U}$.

Now we prove that the theorem holds for $\mathfrak{F}$.

Since $H/H\in \mathfrak{U}$, by the assertion proved above and Lemma
2.2, we see that $H$ is supersoluble. In particular, $H$ is soluble
and hence $F^{*}(H)=F(H)$. Now by using Lemma 2.9, we obtain that
$G\in \mathfrak{F}$. This completes the proof of the theorem.

\section{New Characterization of $p$-nilpotent groups}

\begin{lemma} Let $G$ be a group and $p$ a prime divisor of
$|G|$ with $(|G|,(p-1)(p^{2}-1)\cdot \cdot \cdot (p^{n}-1))=1$ for
some integer $n\geq 1$. Suppose $P$ is a Sylow $p$-subgroup of $G$
and every $n$-maximal subgroup of $P$(if exists) has a $p$-nilpotent
supplement in $G$. Then $G$ is $p$-nilpotent.
\end{lemma}

{\noindent\bf Proof.} Assume that $p^{n+1}\mid |G|$. Let $P_{n1}$ be
an $n$-maximal subgroup of $P$. By hypothesis, $P_{n1}$ has a
$p$-nilpotent supplement $T_{1}$ in $G$. Let $K_{1}$ be a normal
Hall $p'$-subgroup of $T_{1}$. Obviously, $K_{1}$ is a Hall
$p'$-subgroup of $G$. Hence $G=P_{n1}T_{1}=P_{n1}N_{G}(K_{1})$. We
claim that $K_{1}\unlhd G$. Indeed, if $K_{1}\ntrianglelefteq G$,
then $N_{P}(K_{1})=N_{G}(K_{1})\cap P\neq P$ since $T_{1}\subseteq
N_{G}(K_{1})$. Therefore, there exists a maximal subgroup $P_{2}$ of
$P$ such that $N_{P}(K_{1})\leq P_{2}$. Let $P_{n2}$ be an
$n$-maximal subgroup of $P$ contained in $P_{2}$. Since $P=P\cap
G=P\cap P_{n1}N_{G}(K_{1})=P_{n1}(P\cap
N_{G}(K_{1}))=P_{n1}N_{P}(K_{1})$, we have $P_{n1}\neq P_{n2}$. By
hypothesis, $P_{n2}$ has a $p$-nilpotent supplement in $G$. With the
same discussion as above, we can find a Hall $p'$-subgroup $K_{2}$
of $G$ such that $G=P_{n2}N_{G}(K_{2})=P_{2}N_{G}(K_{2})$. If $p=2$,
then by Lemma 2.10, $K_{1}$ conjugates with $K_{2}$ in $G$. If
$p>2$, then $G$ is soluble by Feit-Thompson Theorem. Hence, $K_{1}$
also conjugates with $K_{2}$ in $G$. This means that there exists an
element $g\in P_{2}$, such that $(K_{2})^{g}=K_{1}$. Then
$G=(P_{2}N_{G}(K_{2}))^{g}=P_{2}N_{G}(K_{1})$. Hence, $P=P\cap
G=P\cap P_{2}N_{G}(K_{1})=P_{2}(P\cap
N_{G}(K_{1}))=P_{2}N_{P}(K_{1})=P_{2}$. This contradiction shows
that $p^{n+1}\nmid |G|$. Thus $G$ is $p$-nilpotent by Lemma 2.11.

\begin{lemma} Let $G$ be a group and $p$ a prime divisor of $|G|$ with
$(|G|,(p-1)(p^{2}-1)\cdot \cdot \cdot (p^{n}-1))=1$ for some integer
$n\geq 1$. Suppose that $G$ has a Sylow $p$-subgroup $P$ such that
every $n$-maximal subgroup of $P$ (if exists) either has a
$p$-nilpotent supplement or is $\mathfrak{U_{\mathrm
s}}$-quasinormal in $G$, then $G$ is $p$-nilpotent.
\end{lemma}

{\noindent\bf Proof.} Suppose the Lemma is false and let $G$ be a
counterexample of minimal order. By Lemma 2.11, we have $p^{n+1}\mid
|G|$. Hence $P$ has a non-trivial $n$-maximal subgroup. We proceed
via the following steps:

\textsl{(1) $O_{p'}(G)=1$.}

If $O_{p'}(G)\neq 1$. Then we may choose a minimal normal subgroup
$N$ of $G$ such that $N\leq O_{p'}(G)$. Clearly,
$(|G/N|,(p-1)(p^{2}-1)\cdot \cdot \cdot (p^{n}-1))=1$ and $PN/N$ is
a Sylow $p$-subgroup of $G/N$. Assume that $L/N$ is an $n$-maximal
subgroup of $PN/N$. Then, clearly, $L/N=M_{p}N/N$, where $M_{p}$ is
an $n$-maximal subgroup of $P$. By hypothesis, $M_{p}$ either has a
$p$-nilpotent supplement or is $\mathfrak{U_{\mathrm
s}}$-quasinormal in $G$. By Lemma 2.7(1) and Lemma 2.2(3), we see
that $G/N$ (with respect to $PN/N$) satisfies the hypothesis. The
minimal choice of $G$ implies that $G/N$ is $p$-nilpotent and
consequently $G$ is $p$-nilpotent, a contradiction.

\textsl{(2) $P$ has a maximal subgroup $P_{1}$ such that $P_{1}$ has
no $p$-nilpotent supplement in $G$} (This follows from Lemma 4.1).

\textsl{(3) $G$ is soluble.}

Suppose that $G$ is not soluble. Then $p=2$ by the well known
Feit-Thompson Theorem. Assume that  $O_2(G)\neq 1$. By Lemma 2.7 and
Lemma 2.2(2), $G/O_2(G)$ satisfies the hypothesis. Hence $G/O_2(G)$
is $2$-nilpotent. It follow that $G$ is soluble, a contradiction.
Now assume that $O_2(G)=1$. Then $(P_n)_G=1$, where $P_n$ is an
$n$-maximal subgroup of $P$. Since $P_{n}$ has no $p$-nilpotent
supplement in $G$, $P_n$ is $\mathfrak{U}_{\mathrm s}$-quasinormal
in $G$ by the hypothesis. Hence there exists $K\unlhd G$ such that
$P_nK$ is $s$-permutable in $G$ and $P_n\cap K\leq
Z^\mathfrak{U}_\infty (G)$. If $K=1$, then $P_n\leq O_2(G)=1$ by
Lemma 2.3(1), a contradiction. Thus, $K\neq 1$. If $
Z^\mathfrak{U}_\infty (G)\neq 1$, then there exists a minimal normal
subgroup $H$ of $G$ contained in $Z^\mathfrak{U}_\infty (G)$. Hence
$H$ is of prime power order. This is impossible since $O_{2'}(G)=1$
and $O_2(G)=1$.  Hence $P_n\cap K=1$ and so $2^{n+1}\nmid |K|$. Then
by Lemma 2.11, $K$ has a normal Hall $2'$-subgroup $T$. Since $T$
char $K\unlhd G$, $T\unlhd G$. It follows from (1) that $T=1$.
Consequently, $K\leq O_2(G)=1$, a contradiction again. Hence (3)
holds.

\textsl{(4) $N=O_p(G)$ is the only minimal normal subgroup of $G$
and $G=[N]M$, where $M$ is a maximal subgroup of $G$ and $M$ is
$p$-nilpotent. }

 Let $N$ be a minimal normal subgroup of $G$. By (1)
and (3), $N$ is an elementary abelian $p$-group and $N\leq O_p(G)$.
By Lemma 2.7(1) and Lemma 2.2(2), $G/N$ satisfies the hypothesis and
so $G/N$ is $p$-nilpotent. Since the class of all $p$-nilpotent
groups is a saturated formation, $N$ is the unique minimal normal
subgroup of $G$ and $\Phi (G)=1$. Hence $O_p(G)=N=C_{G}(N)$, and
consequently $G=[N]M$, where $M$ is a $p$-nilpotent maximal subgroup
of $G$. Thus (4) holds.

\textsl{(5) The final contradiction.}

Let $P_{n}$ be an $n$-maximal subgroup of $P$ such that $P_n\leq P_1
$. Then $P_{n}$ has also no $p$-nilpotent supplement in $G$. Hence
there exists a normal subgroup $K$ of $G$ such that $P_{n}K$ is
$s$-permutable in $G$ and $(P_{n}\cap K)(P_{n})_{G}/(P_{n})_{G}\leq
Z_\infty ^\frak{F}(G/(P_{n})_{G})$. We claim that $(P_{n})_{G}=1$.
Indeed, if $(P_{n})_{G}\neq 1$, then by (2),
$O_{p}(G)=N=(P_{n})_{G}$. Hence $G=NM=(P_{n})_{G}M=P_{n}M$, which
contradicts (2). Therefore, $P_{n}\cap K\leq Z_\infty ^\frak{U}(G)$.
If $K=1$, then $P_{n}$ is $s$-permutable in $G$, and so $P_{n}\leq
O_{p}(G)=N$ and $O^{p}(G)\leq N_{G}(P_{n})$ by Lemma 2.3. Hence
$1\neq P_{n}\leq P_{n}^{G}=P_{n}^{O^{p}(G)P}=P_{n}^{P}=(P_{n}\cap
N)^{P}\leq (P_{1}\cap N)^{P}=P_{1}\cap N \leq N$. On the other hand,
obviously, $N\leq P_{n}^{G}$. Thus $N=P_{n}^{G}=P_{1}\cap N$. It
follows that $N\leq P_{1}$, and so $G=NM=P_{1}M$. This means that
$P_{1}$ has a $p$-nilpotent supplement in $G$. This contradiction
shows that $K\neq 1$. If $P_{n}\cap K=1$, then $p^{n+1}\nmid |K|$.
By Lemma 2.11, $K$ is $p$-nilpotent and $K_{p'}\leq O_{p'}(G)=1$ by
(1). Hence $K=N=O_{p}(G)$. It follows from Lemma 2.3(1) that
$P_{n}K=K$ and so $P_{n}\cap K\neq 1$, a contradiction. Hence
$P_{n}\cap K\neq 1$. This means that $Z_\infty ^\frak{U}(G)\neq 1$
and so $N\leq Z_\infty ^\frak{U}(G)$. Consequently,
$|N|=|O_p(G)|=p$. Therefore, $G/N\cong G/C_{G}(N)$ is isomorphic
with some subgroup of $Aut(N)$ of order $p-1$. Since
$(|G|,(p-1)(p^{2}-1)\cdot \cdot \cdot (p^{n}-1))=1$, $G/N=1$.
Consequently, $G=N$ is an elementary abelian $p$-group. The final
contradiction completes the proof.

\begin{theorem} \textbf{Let $p$ be a prime and $G$ a group. Suppose that
$(|G|,(p-1)(p^{2}-1)\cdot \cdot \cdot (p^{n}-1))=1$ for some integer
$n\geq 1$. Then $G$ is $p$-nilpotent if and only if $G$ has a
normal subgroup $E$ such that $G/E$ is $p$-nilpotent and every
$n$-maximal subgroup of $P$(if exists) either has a $p$-nilpotent
supplement or is $\mathfrak{U_{\mathrm s}}$-quasinormal in $G$,
where $P$ is a Sylow $p$-subgroup of $E$.}
\end{theorem}
{\noindent\bf Proof.} The necessity is obvious. We only need to
prove the sufficiency. Suppose it is false and let $G$ be a
counterexample of minimal order. By Lemma 2.7(2) and Lemma 2.2(4),
every $n$-maximal subgroup of $P$ either has a $p$-nilpotent
supplement or is $\mathfrak{U_{\mathrm s}}$-quasinormal in $E$.
Hence $E$ is $p$-nilpotent by Lemma 4.2. Then, $E\neq G$. Let $T$ be
a normal Hall $p'$-subgroup of $E$. Clearly, $T\unlhd G$. We proceed
the proof via the following steps:

\textsl{(1) $T=1$, and so $P=E\unlhd G$.}

Suppose that $T\neq 1$. Since $T$ is a normal Hall $p'$-subgroup of
$E$ and $E\unlhd G$, then $T\unlhd G$. We show that $G/T$ (with
respect to $E/T$) satisfies the hypothesis. Indeed,
$(G/T)/(E/T)\simeq G/E$ is $p$-nilpotent and $E/T=PT/T$ is a
$p$-group. Suppose that $M_{n}/T$ is an $n$-maximal subgroup of
$PT/T$ and $P_{n}=M_{n}\cap P$. Then $P_{n}$ is an $n$-maximal
subgroup of $P$ and $M_{n}=P_{n}T$. By the hypothesis, $P_{n}$
either has a $p$-nilpotent supplement or is $\mathfrak{U_{\mathrm
s}}$-quasinormal in $G$. By Lemma 2.7(1) and Lemma 2.2(3),
$M_{n}/T=P_{n}T/T$ either has a $p$-nilpotent supplement or is
$\mathfrak{U_{\mathrm s}}$-quasinormal in $G/T$. The minimal choice
of $G$ implies that $G/T$ is $p$-nilpotent. This implies that $G$ is
$p$-nilpotent. This contradiction shows $T=1$. Hence $P=E\unlhd G$.

\textsl{(2) Let $Q$ be a Sylow $q$-subgroup of $G$, where $q$ is a
prime divisor of $|G|$ with $q\neq p$. Then $PQ=P\times Q$.}

By (1), $P=E\unlhd G$, $PQ$ is a subgroup of $G$. By Lemma 2.7(2)
and Lemma 2.2(4), every $n$-maximal subgroup of $P$ either has a
$p$-nilpotent supplement or is $\mathfrak{U_{\mathrm
s}}$-quasinormal in $PQ$. By using Lemma 4.2, we have that $PQ$ is
$p$-nilpotent. Hence $Q\unlhd PQ$ and thereby $PQ=P\times Q$.

\textsl{(3) The final contradiction.}

From (2), we have $O^p(G)\leq C_G(P)$. \textbf{This induces that $E=P\leq
Z_\infty(G)$. Therefore $G$ is $p$-nilpotent.} The final contradiction completes the proof.

\begin{theorem} Let $G$ be a finite group and $p$ a prime divisor of
$|G|$ with $(|G|,p-1)=1$. Then $G$ is $p$-nilpotent if and only if
 $G$ has a soluble normal subgroup $H$ of $G$ such that
$G/H$ is $p$-nilpotent and every maximal subgroup of every Sylow
subgroup of $F(H)$ is $\mathfrak{U_{\mathrm s}}$-quasinormal in $G$.
\end{theorem}

{\noindent\bf Proof.} The necessity is obvious. We only need to
prove the sufficiency.  Suppose that it is false and let $G$ be a
counterexample with $|G||H|$ is  minimal. Let $P$ be an arbitrary
given Sylow $p$-subgroup of $F(H)$. Clearly, $P\unlhd G$. We proceed
the proof as follows.

\textsl{(1) $\Phi(G)\cap P=1$.}

If not, then $1\neq \Phi(G)\cap P \unlhd G$. Let $R=\Phi(G)\cap P$.
Clearly, $(G/R)/(H/R)\simeq G/H\in \mathfrak{F}$. By
Gasch$\ddot{u}$tz theorem (see [11, III, Theorem 3.5]), we have that
$F(H/R)=F(H)/R$. Assume that $P/R$ is a Sylow $\textit{p}$-subgroup
of $F(H/R)$ and $P_{1}/R$ is a maximal subgroup of $P/R$. Then $P$
is a Sylow $p$-subgroup of $F(G)$ and $P_{1}$ is a maximal subgroup
of $P$. By Lemma 2.2(2) and the hypothesis, $P_{1}/R$ is
$\mathfrak{U_{\mathrm s}}$-quasinormal in $G/R$. Now, let $Q/R$ be a
maximal subgroup of some Sylow $\textit{q}$-subgroup of
$F(H/R)=F(H)/R$, where $q\neq p$. Then $Q=Q_{1}R$, where $Q_{1}$ is
a maximal subgroup of the Sylow $\textit{q}$-subgroup of $F(H)$. By
hypothesis, $Q_{1}$ is $\mathfrak{U_{\mathrm s}}$-quasinormal in
$G$. Hence $Q/R=Q_{1}R/R$ is $\mathfrak{U_{\mathrm s}}$-quasinormal
in $G/R$ by Lemma 2.2(3). This shows that $(G/R, H/R)$ satisfies the
hypothesis. The minimal choice of $(G,H)$ implies that $G/R$ is
$p$-nilpotent. It follows that $G$ is $p$-nilpotent, a
contradiction. Hence (1) holds.

\textsl{(2) $P=\langle x_{1}\rangle \times \langle x_{2}\rangle
\times$ $\cdot \cdot \cdot \times \langle x_{m}\rangle$, where every
$\langle x_{i}\rangle$ $(i\in \{1\cdot \cdot \cdot m\})$ is a normal
subgroup of $G$ with order $p$.}

By (1) and Lemma 2.5, $P=R_{1}\times R_{2}\times$ $\cdot \cdot \cdot
\times R_{m}$, where $R_{i} \ (i\in \{1\cdot \cdot \cdot m\})$ is a
minimal normal subgroup of $G$. We now prove that $R_{i}$ is of
order $\textit{p}$, for $i\in \{1\cdot \cdot \cdot m\}$.

Assume that $|R_{i}|>\textit{p}$, for some $\textit{i}$. Without
loss of generality, we let $|R_{1}|>\textit{p}$ and $R^{*}_{1}$ be a
maximal subgroup of $R_{1}$. Then, $R^{*}_{1}\neq 1$ and
$R^{*}_{1}\times R_{2}\times$ $\cdot \cdot \cdot \times R_{m}=P_{1}$
is a maximal subgroup of $P$. Put $T=R_{2}\times$ $\cdot \cdot \cdot
\times   R_{m}$. Then, clearly, $(P_{1})_{G}=T$. By hypothesis,
$P_{1}$ is $\mathfrak{U_{\mathrm s}}$-quasinormal in $G$. Hence by
Lemma 2.2(1), there exists a normal subgroup $N$ of $G$ such that
$(P_{1})_{G}\leq N$, $P_{1}N$ is $s$-permutable in $G$ and
$P_{1}/(P_{1})_{G} \cap N/(P_{1})_{G}\leq Z_\infty ^\frak{U}
(G/(P_{1})_{G})$. Assume that $P_{1}/(P_{1})_{G} \cap
N/(P_{1})_{G}\neq 1$. Let $Z_\infty ^\frak{U}
(G/(P_{1})_{G})=V/(P_{1})_{G}=V/T$. Then $P_{1}\cap N\leq V$ and
$P/T \cap V/T \unlhd G/T$. Since $P\cap V\geq P_{1}\cap N\cap V\geq
P_{1}\cap N> (P_{1})_{G}=T$, $P/T \cap V/T\neq 1$. As $P/T\simeq
R_{1}$ and $R_{1}$ is a minimal normal subgroup of $G$, we have
$P/T\subseteq V/T$. This implies that $|R_{1}|=|P/T|=p$. This
contradiction shows that $P_{1}\cap N=(P_{1})_{G}=T$. Consequently,
$P_{1}N=R^{*}_{1}TN=R^{*}_{1}N$ and $R^{*}_{1}\cap N=1$. Since
$R_{1}\cap N\unlhd G$, $R_{1}\cap N=1$ or $R_{1}\cap N=R_{1}$. If
$R_1\cap N=R_1$, then $R_1^*\subseteq R_1\subseteq N$, which
contradicts $R^{*}_{1}\cap N=1$. Hence $R_{1}\cap N=1$. It follows
that $R^{*}_{1}=R^{*}_{1}(R_{1}\cap N)=R_{1}\cap R^{*}_{1}N$ is
$s$-permutable in $G$. Thus $O^{p}(G)\leq N_{G}(R_{1}^{*})$ by Lemma
2.3(2). This induces that for every maximal subgroup $R_{1}^{*}$ of
$R_{1}$, we have that $|G:N_{G}(R_{1}^{*})|=p^{\alpha}$, where
$\alpha$ is an integer. Let $\{ R_{1}^{*},R_{2}^{*},\cdot \cdot
\cdot ,R_{t}^{*}\}$ be the set of all maximal subgroups of $R_{1}$.
Then $p$ divides $t$. This contradicts to [11, III, Theorem 8.5(d)].
Thus (2) holds.

\textsl{(3) $G/F(H)$ is $p$-nilpotent.}

By (2), $F(H)=\langle y_{1}\rangle \times \langle y_{2}\rangle
\times$ $\cdot \cdot \cdot \times \langle y_{n}\rangle$, where
$\langle y_{i}\rangle \ (i\in \{1\cdot \cdot \cdot n\})$ is a normal
subgroup of $G$ of order $p$. Since $G/C_{G}(\langle y_{i}\rangle)$
is isomorphic with some subgroup of $Aut(\langle y_{i}\rangle)$,
$G/C_{G}(\langle y_{i}\rangle)$ is cyclic. Hence, $G/C_{G}(\langle
y_{i}\rangle)$ is $p$-nilpotent for every $i$. It follows that
$G/\cap^{n}_{i=1}C_{G}(\langle y_{i}\rangle)$ is $p$-nilpotent.
Obviously, $C_{G}(F(G))=\cap^{n}_{i=1}C_{G}(\langle y_{i}\rangle)$.
Hence $G/C_{G}(F(G))$ is $p$-nilpotent. Consequently, $G/(H\cap
C_{G}(F(G)))=G/C_{H}(F(H))$ is $p$-nilpotent. Since $F(H)$ is
abelian, $F(H)\leq C_{H}(F(H))$. On the other hand, $C_{H}(F(H))\leq
F(H)$ since $H$ is soluble. Thus $F(H)=C_{H}(F(H))$ and so $G/F(H)$
is $p$-nilpotent.

\textsl{(4) If $K$ is a minimal normal subgroup of $G$ contained in
$H$, then $K\subseteq F(H)$ and  $G/K$ is $p$-nilpotent.}

Let $K$ be an arbitrary minimal normal subgroup of $G$ contained in
$H$. Then $K$ is an elementary abelian $\textit{p}$-group for some
prime $\textit{p}$ since $H$ is soluble. Hence $K\leq F(H)$. By
Lemma 2.2(2) and (3), we see that $G/K$ (with respect to $H/K$)
satisfies the hypothesis. The minimal choice of $(G,H)$ implies that
$G/K$ is $p$-nilpotent.

\textsl{(5) The final contradiction.}

Since the class of all $p$-nilpotent groups is a saturated
formation, by (2) and (4), we see that $K=F(H)=\langle x\rangle$ is
the unique minimal normal subgroup of $G$ contained in $H$, where
$\langle x\rangle$ is a cyclic group of order $p$ for some prime
$p$. Since $G/K$ is $p$-nilpotent, it has a normal $p$-complement
$L/K$. By Schur-Zassenhaus Theorem, $L=G_{p'}K$, where $G_{p'}$ is a
Hall $p'$-subgroup of $G$. Since $p$ is the prime divisor of $|G|$
with $(|G|,p-1)=1$ and $N_{L}(K)/C_{L}(K)\simeq Aut(K)$ is a
subgroup of a cyclic group of order $p-1$, we see that
$N_{L}(K)=C_{L}(K)$. Then, by Burnside Theorem (see [14, (10.1.8)]),
we have that $L$ is $p$-nilpotent. Then $G_{p'}$ $\emph{char}$
$L\unlhd G$, so $G_{p'}\unlhd G$. Hence $G$ is $p$-nilpotent. The
final contradiction completes the proof.

\

\noindent{\bf Acknowledgements.}  The authors are very grateful to
the helpful suggestions of the referee.


\begin{thebibliography}{1}
\bibitem{} M. Assad, Piroska Cs\"org\H{o}, Characterization of
Finite Groups With Some $S$-quasinormal subgroups, Monatsh. Math.
\textbf{146} (2005), 263-266.
\bibitem{} A. Ballester-Bolinches, M.C.Pedraza-Aguilera, Sufficient
conditions for supersolubility of finite groups, J. Pure Appl.
Algebra, \textbf{127} (1998), 113-118.
\bibitem{} F. Gross, Conjugacy of odd order Hall subgroups, Bull. London, math. soc, \textbf{19} (1987), 311-319.
\bibitem{} W. Guo, The Theory of Classes of Groups, Science
Press/Kluwer, Academic Publishers, Beijing New York, 2000.
\bibitem{} W. Guo, On $\mathfrak{F}$-supplemented subgroups of
finite group, Manuscripta Math. {\bf 127} (2008), 139-150.
\bibitem{} W. Guo, K. P. Shum and F. Xie, Finite groups with some weakly
$S$-supplemented subgroups, Glasgow Math. J. \textbf{53} (2011),
211-222.
\bibitem{} W. Guo, K. P. Shum and A. N. Skiba, $X$-semipermutable subgroup of finite groups. J.
Algebra, \textbf{315} (2007), 31-41.
\bibitem{} W. Guo, A. N. Skiba, Finite groups with given $s$-embedded
and $n$-embedded subgroups, J.Algebra, \textbf{321} (2009),
2843-2860.
\bibitem{} W. Guo, A. N. Skiba, New criterions of existence and
conjugacy of Hall subgroups of finite groups. Proceedings of AMS,
\textbf{139} (2011), 2327-2336.
\bibitem{} J. Huang, On $\mathfrak{F_{\mathrm s}}$-quasinormal Subgroups of Finite Groups,
Comm. Algebra, \textbf{38} (2010), 4063-4076.
\bibitem{} B. Huppert, Endliche Gruppen I, Springer-Verlag, Berlin,
1967.
\bibitem{} B. Huppert, N. Blackburn, Finite Groups III, Springer-Verlag, Berlin,
New York, 1982.
\bibitem{} O. Kegel, Sylow-Gruppen and Subnormalteiler endlicher
Gruppen, Math. Z., \textbf{78} (1962), 205-221.
\bibitem{} D. J. S. Robinson, A course in the Theory of
Groups. New York: Springer, 1982.
\bibitem{} L. A. Shemetkov, Formations of Finite groups Moscow, Nauka,
1978.
\bibitem{} L.A. Shemetkov, A.N. Skiba, Formations of Algebraic
Systems, Moscow, Nauka, 1989.
\bibitem{} A.N. Skiba, On weakly $s$-permutable subgroups of finite
groups, J. Algebra, \textbf{315} (2007), 192-209.
\bibitem{} H. Wielandt, Subnormal subgroups and permutation groups,
Lectures given at the Ohio State University, Columbus, Ohio, 1971.
\end{thebibliography}
\end{document}